\documentclass[12pt]{amsart}
\usepackage{amssymb, hyperref}
\usepackage[margin=0.8in]{geometry}

\newcommand{\nc}{\newcommand}
\nc{\set}[2]{\{\, #1 : #2\,\}}
\nc{\st}{\op{star}}
\nc{\tPa}[6]{\bibitem{#1} {#2}, \emph{#3}, {#4}, to appear (#5 pages).}
\nc{\sPa}[5]{\bibitem{#1} {#2}, \emph{#3}, {#4}, submitted.}
\nc{\FinSeqs}[1]{{{#1}^{<\alephes}}}
\nc{\rest}[1]{\mid_{#1}}
\nc{\Bgp}{{\Z^\N}}
\nc{\Arh}{Arhangel'ski\u{\i}}
\nc{\grbl}{{\mbox{\textit{\tiny gp}}}}
\long\def\forget#1\forgotten{}
\nc{\issuenumber}{38}
\nc{\issuemonth}{Nov}
\nc{\issueyear}{2015}

\nc{\Ga}{\Gamma}
\nc{\sub}{\subseteq}
\nc{\R}{\mathbb{R}}
\nc{\Cp}{\mathrm{C}_\mathrm{p}}
\nc{\Op}{\mathrm{O}}
\nc{\Fin}{{[\N]^{<\aleph_0}}}
\nc{\alephes}{{\aleph_0}}
\nc{\ed}{
\forget

\section{Unsolved problems from earlier issues}

\begin{issue}Is $\binom{\Omega}{\Gamma}=\binom{\Omega}{\Tau}$?\end{issue}%\stepcounter{issue}
\begin{issue}Is $\ufin(\Op,\Omega)=\sfin(\Gamma,\Omega)$?And if not, does $\ufin(\Op,\Gamma)$ imply
$\sfin(\Gamma,\Omega)$?\end{issue}%\stepcounter{issue}
\stepcounter{issue}
\begin{issue}Does $\sone(\Omega,\Tau)$ imply $\ufin(\Gamma,\Gamma)$?\end{issue}
\stepcounter{issue}
\begin{issue}Is there, in ZFC, an uncountable set satisfying $\sfin(\cB,\cB)$?\end{issue}
\stepcounter{issue}
\begin{issue}Does $X \nin \NON(\cM)$ and $Y\nin\mathsf{D}$ imply that $X\cup Y\nin \COF(\cM)$?\end{issue}
\begin{issue}[CH]Is $\split(\Lambda,\Lambda)$ preserved under finite unions?\end{issue}
\begin{issue}Is $\cov(\cM)=\fo$? (See the definition of $\fo$ in that issue.)\end{issue}
\stepcounter{issue}
\begin{issue}Could there be a Baire metric space $M$ of weight $\aleph_1$ and a partition
$\mathcal{U}$ of $M$ into $\aleph_1$ meager sets where for each ${\mathcal U}'\subset\mathcal U$,
$\bigcup {\mathcal U}'$ has the Baire property in $M$?\end{issue}
\stepcounter{issue} %% no problem in Issue 13
\begin{issue}Is there, in ZFC, a set of reals $X$ of cardinality $\fd$ such that all
finite powers of $X$ have Menger's property $\sfin(\Op,\Op)$?\end{issue}
\stepcounter{issue}
\begin{issue}[MA]Is there a set $X\sbst\bbR$, of cardinality continuum, satisfying $\sone(\BO,\BG)$?\end{issue}
\begin{issue}[CH]Is there a totally imperfect $X$ satisfying $\ufin(\Op,\Gamma)$
that can be mapped continuously onto $\Cantor$?\end{issue}
\begin{issue}[CH]Is there a Hurewicz $X$ such that $X^2$ is Menger but not Hurewicz?\end{issue}
\begin{issue}Does the Pytkeev property of $C_p(X)$ imply that $X$ has Menger's property?\end{issue}
\begin{issue}Does every hereditarily Hurewicz space satisfy $\sone(\BG,\BG)$?\end{issue}
\begin{issue}[CH]Is there a Rothberger-bounded $G\le\Bgp$ such that $G^2$ is not Menger-bounded?\end{issue}
\begin{issue}Let $\cW$ be the van der Waerden ideal. Are $\cW$-ultrafilters closed under products?\end{issue}
\begin{issue}Is the $\delta$-property equivalent to the $\gamma$-property $\binom{\Omega}{\Gamma}$?\end{issue}
\stepcounter{issue}\stepcounter{issue}

\forgotten

\general\end{document}}

\nc{\fbx}[1]{\fbox{$#1$}}\nc{\nop}{$\times$}\nc{\fbn}{\!\!\fbox{\!\nop\!}\!\!}
\nc{\yup}{\checkmark}\nc{\fby}{\!\!\fbox{\!\yup\!}\!\!}\nc{\mbq}{\mb{?}}
\nc{\mb}[1]{{\mbox{\textbf{#1}}}}\nc{\smb}[1]{{\!\!\mb{#1}\!\!}}\nc{\x}{\times}
\nc{\<}{\left <}\nc{\Cantor}{{\{0,1\}^\N}}\nc{\oo}{\infty}
\nc{\NR}{{\bbR^\N}}\nc{\Iff}{\Leftrightarrow}\nc{\mypar}[1]{\par\medskip\noindent\textbf{#1.}}
\nc{\roth}{{[\N]^{\alephes}}}\nc{\sr}[2]{{\txt{$#1$\\$#2$}}}
\nc{\gpbl}{{\mbox{\textit{\tiny gp}}}}\nc{\fx}{\mathfrak{x}}\nc{\fb}{\mathfrak{b}}
\nc{\fg}{\mathfrak{g}}\nc{\fc}{\mathfrak{c}}\nc{\fd}{\mathfrak{d}}
\nc{\fp}{\mathfrak{p}}\nc{\fs}{\mathfrak{s}}\nc{\ADD}{{\mathsf   {ADD}}}
\nc{\COV}{{\mathsf   {COV}}}\nc{\NON}{{\mathsf   {NON}}}\nc{\COF}{{\mathsf   {COF}}}
\nc{\sseq}[1]{\{#1 : n\in\N\}}\nc{\Impl}{\Rightarrow}
\nc{\upannouncement}[1]{[\S\ref{#1} above]}\nc{\dnannouncement}[1]{[\S\ref{#1} below]}
\nc{\E}{\exists}\nc{\cI}{\mathcal{I}}\nc{\cN}{\mathcal{N}}\nc{\cP}{\mathcal{P}}
\nc{\cA}{\mathcal{A}}\nc{\cM}{\mathcal{M}}\nc{\Null}{\mathcal{N}}
\nc{\op}{\operatorname}\nc{\cov}{\mathsf{cov}}\nc{\add}{\mathsf{add}}
\nc{\cof}{\mathsf{cof}}\nc{\cf}{\mathsf{cf}}\nc{\non}{\mathsf{non}}\nc{\spst}{\supseteq}
\nc{\CH}{the Continuum Hypothesis}\nc{\bbR}{\mathbb{R}}\nc{\Q}{\mathbb{Q}}
\nc{\EdNote}[1]{\par\medskip\noindent\textbf{#1.}}\nc{\fo}{\mathfrak{od}}
\nc{\cl}[1]{\overline{#1}}\nc{\impl}{\rightarrow}\nc{\arrays}{{{\{0,1\}}^{\N\x\N}}}
\nc{\w}{\omega}\nc{\ft}{\mathfrak{t}}\nc{\h}{\mathfrak{h}}\nc{\Cite}[2]{{\cite[#1]{#2}}}
\renewcommand{\split}{\mathsf{Split}}\nc{\bq}{\begin{quote}}\nc{\eq}{\end{quote}}
\nc{\cK}{\mathcal{K}}\nc{\cB}{\mathcal{B}}\nc{\BG}{\cB_\Gamma}
\nc{\BL}{\cB_\Lambda}\nc{\BT}{\cB_\Tau}\nc{\BTstar}{\cB_{\Tau^*}}\nc{\BO}{\cB_\Omega}
\nc{\CG}{C_\Gamma}\nc{\CL}{C_\Lambda}\nc{\CT}{C_\Tau}\nc{\CTstar}{C_{\Tau^*}}
\nc{\CO}{C_\Omega}\nc{\sone}{\mathsf{S}_1}\nc{\sfin}{\mathsf{S}_\mathrm{fin}}
\nc{\Sc}{\mathsf{S}_c}\nc{\ufin}{\mathsf{U}_\mathrm{fin}}\nc{\gone}{\mathsf{G}_1} \nc{\gfin}{\mathsf{G}_\mathrm{fin}}\nc{\seq}[1]{\{#1\}_{n\in\N}}\nc{\Un}{\bigcup}
\nc{\nin}{\not\in}\nc{\cF}{\mathcal{F}}\nc{\cG}{\mathcal{G}}\nc{\cU}{\mathcal{U}}
\nc{\cV}{\mathcal{V}}\nc{\cW}{\mathcal{W}}\nc{\fU}{\mathfrak{U}}\nc{\fu}{\mathfrak{u}}
\nc{\fV}{\mathfrak{V}}\nc{\fW}{\mathfrak{W}}\nc{\psin}{pseudo-intersection}
\nc{\NN}{{\N^\N}}\nc{\N}{\mathbb{N}}\nc{\bbN}{\mathbb{N}}\nc{\Z}{\mathbb{Z}}
\nc{\as}{\subseteq^*}\nc{\sm}{\setminus}\nc{\sbst}{\subseteq}
\nc{\by}[2]{\par\hfill\emph{#1}, #2}\nc{\nby}[1]{\par\hfill\emph{#1}}\nc{\Tau}{\mathrm{T}}
\nc{\CE}{\textsc{CE}}
\newtheorem{thm}{Theorem}[section]\nc{\bthm}{\begin{thm}} \nc{\ethm}{\end{thm}}
\newtheorem{prop}[thm]{Proposition}\nc{\bprp}{\begin{prop}} \nc{\eprp}{\end{prop}}
\newtheorem{fact}[thm]{Fact}\nc{\bfct}{\begin{fact}} \nc{\efct}{\end{fact}}
\newtheorem{prob}[thm]{Problem}\nc{\bprb}{\begin{prob}} \nc{\eprb}{\end{prob}}
\newtheorem{lem}[thm]{Lemma}\nc{\blem}{\begin{lem}} \nc{\elem}{\end{lem}}
\newtheorem{claim}[thm]{Claim}\nc{\bclm}{\begin{claim}} \nc{\eclm}{\end{claim}}
\newtheorem{cor}[thm]{Corollary}\nc{\bcor}{\begin{cor}} \nc{\ecor}{\end{cor}}
\newtheorem{conj}[thm]{Conjecture}\nc{\bcnj}{\begin{conj}} \nc{\ecnj}{\end{conj}}
\theoremstyle{definition}\newtheorem{defn}[thm]{Definition}\nc{\bdfn}{\begin{defn}} \nc{\edfn}{\end{defn}}
\theoremstyle{remark}\newtheorem{rem}[thm]{Remark}\nc{\brem}{\begin{rem}} \nc{\erem}{\end{rem}}
\newtheorem{cnv}[thm]{Convention}\nc{\bcnv}{\begin{cnv}} \nc{\ecnv}{\end{cnv}}
\newtheorem{exam}[thm]{Example}\nc{\bexm}{\begin{exam}} \nc{\eexm}{\end{exam}}
\newtheorem{issue}{Issue}\nc{\bpf}{\begin{proof}} \nc{\epf}{\end{proof}}
\nc{\be}{\begin{enumerate}}\nc{\ee}{\end{enumerate}}\nc{\bi}{\begin{itemize}}
\nc{\ei}{\end{itemize}}\nc{\itm}{\item}
\nc{\general}{\small\vfill\par\noindent\hrulefill\par
\noindent\textbf{Previous issues.} 
\url{http://front.math.ucdavis.edu/search?\&t=\%22SPM+Bulletin\%22}
\\[0.1cm]
\textbf{Contributions and free subscription.} Email \url{tsaban@math.biu.ac.il}.
}

\nc{\link}[1]{\par\hfill{\url{#1}}}
\nc{\fillin}{{\Huge To be completed}}
\nc{\arXivl}[4]{\subsection{#2}{#4}\par\hfill{\arx{#1}}\par\hfill\emph{#3}}
\nc{\arXiv}[3]{\subsection{#2}\mbox{}\par\hfill{\arx{#1}}\par\hfill\emph{#3}}
\nc{\DOIpaper}[5]{\subsection{#2}{#4}\par\hfill{\url{http://dx.doi.org/#1}}\par\hfill\emph{#3}}
\nc{\AMSPaper}[5]{\subsection{#3}{#5}\par\hfill{\url{#1}}\par\hfill\emph{#4}\par\hfill{#2}}
\nc{\nAMSPaper}[4]{\subsection{#2}{#4}\par\hfill{\url{#1}}\par\hfill\emph{#3}}
\nc{\AMS}[3]{\subsection{#1}\mbox{}\par\hfill{\url{#3}}\par\hfill\emph{#2}}
\nc{\SPMBul}{\textbf{$\mathcal{SPM}$ Bulletin}}
\nc{\BulEnd}{\par\bigskip\noindent
Boaz Tsaban\\
\emph{E-mail}: tsaban@math.biu.ac.il\\
\emph{URL}: http://www.cs.biu.ac.il/\~{}tsaban}

\nc{\arx}[1]{\url{http://arxiv.org/abs/#1}}
\nc{\probissue}{\emph{Problem of the issue}}

\title[$\mathcal{SPM}$ Bulletin \textbf{\issuenumber} (\issuemonth{} \issueyear)]{%
$\mathcal{SPM}$ Bulletin\\[0.5cm]
Issue number \issuenumber: \issuemonth{} \issueyear{} \CE{}}

\begin{document}
\maketitle

%\tableofcontents

%\forget
\section{Editor's note}

The \emph{International Conference on Topology} took place in Messina, September 7--11, 2015. 
A substantial portion of the lectures dealt with selection principles. The program and book of abstracts
are available on the conference webpage, \url{http://mat521.unime.it/ictm2015}\,.

A special issue of the journal \emph{Topology and its Applications} will be dedicated to the conference's themes and,
in particular, to selection principles. The guest editors for this issue are Maddalena Bonanzinga and Boaz Tsaban.
Papers meeting the journal's high standards may be submitted, by the end of January (tentative),
to the special issue's email address \url{ictmessina2015@gmail.com}\,. 
The papers will be fully refereed according to the journal's standards.
Attendance in the conference is not a prerequisite for submission; the sole criteria are quality and relevance (in this order).

\medskip

With best regards,

\by{Boaz Tsaban}{tsaban@math.biu.ac.il}

\hfill \texttt{http://www.cs.biu.ac.il/\~{}tsaban}
%\forgotten

\section{Long announcements}

\arXivl{1502.00177}
{Selective versions of chain condition-type properties}
{Leandro Aurichi, Santi Spadaro, Lyubomyr Zdomskyy}
{We study selective and game-theoretic versions of properties like the ccc,
weak Lindel\"ofness and separability, giving various characterizations of them
and exploring connections between these properties and some classical cardinal
invariants of the continuum.}

\arXivl{1502.00178}
{On topological properties of the weak topology of a Banach space}
{Saak Gabriyelyan, Jerzy Kakol, Lyubomyr Zdomskyy}
{Being motivated by the famous Kaplansky theorem we study various sequential
properties of a Banach space $E$ and its closed unit ball $B$, both endowed
with the weak topology of $E$. We show that $B$ has the Pytkeev property if and
only if $E$ in the norm topology contains no isomorphic copy of $\ell_1$, while
$E$ has the Pytkeev property if and only if it is finite-dimensional. We extend
Schl\"uchtermann and Wheeler's result by showing that $B$ is a (separable)
metrizable space if and only if it has countable $cs^\ast$-character and is a
$k$-space. As a corollary we obtain that $B$ is Polish if and only if it has
countable $cs^\ast$-character and is \v{C}ech-complete, that supplements a
result of Edgar and Wheeler.}

\arXivl{1503.04229}
{On metrizable $X$ with $C_p(X)$ not homeomorphic to $C_p(X)\times
  C_p(X)$}
{Miko{\l}aj Krupski and Witold Marciszewski}
{We give two examples of infinite metrizable spaces $X$ such that the space
$C_p(X)$, of continuous real-valued function on $X$ endowed with the pointwise
topology, is not homeomorphic to its own square $C_p(X)\times C_p(X)$. The
first of them is a one-dimensional continuum; the second one is a
zero-dimensional subspace of the real line. Our result answers a long-standing
open question in the theory of function spaces posed by A.V. Arhangel'skii.}

\arXivl{1503.04383}
{On complete metrizability of the Hausdorff metric topology}
{Laszlo Zsilinszky}
{There exists a completely metrizable bounded metrizable space $X$ with
compatible metrics $d,d'$ so that the hyperspace $CL(X)$ of nonempty closed
subsets of $X$ endowed with the Hausdorff metric $H_d$, $H_{d'}$, resp. is
$\alpha$-favorable, $\beta$-favorable, resp. in the strong Choquet game. In
particular, there exists a completely metrizable bounded metric space $(X,d)$
such that $(CL(X),H_d)$ is not completely metrizable.}

\arXivl{1503.06092}
{$G_\delta$ semifilters and $\omega^*$}
{Will Brian and Jonathan Verner}
{The ultrafilters on the partial order $([\omega]^{\omega},\subseteq^*)$ are
the free ultrafilters on $\omega$, which constitute the space $\omega^*$, the
Stone-\v{C}ech remainder of $\omega$. If $U$ is an upperset of this partial
order (i.e., a \emph{semifilter}), then the ultrafilters on $U$ correspond to
closed subsets of $\omega^*$ via Stone duality.
  If, in addition, $U$ is sufficiently "simple" (more precisely, $G_\delta$ as
a subset of $2^\omega$), we show that $U$ is similar to $[\omega]^{\omega}$ in
several ways. First, $\mathfrak{p}_U = \mathfrak{t}_U = \mathfrak{p}$ (this
extends a result of Malliaris and Shelah). Second, if $\mathfrak{d} =
\mathfrak{c}$ then there are ultrafilters on $U$ that are also $P$-filters
(this extends a result of Ketonen). Third, there are ultrafilters on $U$ that
are weak $P$-filters (this extends a result of Kunen).
  By choosing appropriate $U$, these similarity theorems find applications in
dynamics, algebra, and combinatorics. Most notably, we will answer two open
questions of Hindman and Strauss by proving that there is an idempotent of
$\omega^*$ that is both minimal and maximal.}

\arXivl{1503.08467}
{Selective strong screenability and a game}
{Liljana Babinkostova and Marion Scheepers}
{Selective versions of screenability and of strong screenability coincide in a
large class of spaces. We show that the corresponding games are not equivalent
in even such standard metric spaces as the closed unit interval. We identify
sufficient conditions for ONE to have a winning strategy, and necessary
conditions for TWO to have a winning strategy in the selective strong
screenability game.}

\AMSPaper{www.ams.org/journal-getitem?pii=S0894-0347-2015-00830-X}
{M. Malliaris; S. Shelah}
{Cofinality spectrum theorems in model theory, set theory, and general topology}
{Journal of the American Mathematical Society}
{We connect and solve two long-standing open problems in quite different areas: the model-theoretic question of whether $\mathrm{SOP}_2$ is maximal in Keisler's order, and the question from general topology/set theory of whether $\fp=\ft$, the oldest problem on cardinal invariants of the continuum. We do so by showing these problems can be translated into instances of a more fundamental problem which we state and solve completely, using model-theoretic methods.}

\arXivl{1504.01626}
{Menger remainders of topological groups}
{Angelo Bella, Se\c{c}il Tokg\"oz, and Lyubomyr Zdomskyy}
{In this paper we discuss what kind of constrains combinatorial covering
properties of Menger, Scheepers, and Hurewicz impose on remainders of
topological groups. For instance, we show that such a remainder is Hurewicz if
and only it it is $\sigma$-compact. Also, the existence of a Scheepers
non-$\sigma$-compact remainder of a topological group follows from CH and
yields a $P$-point, and hence is independent of ZFC. We also make an attempt to
prove a dichotomy for the Menger property of remainders of topological groups
in the style of Arhangel'skii.}

\arXivl{1505.06238}
{The Whyburn property and the cardinality of topological spaces}
{Santi Spadaro}
{The weak Whyburn property is a generalization of the classical sequential
property that was studied by many authors. A space $X$ is weakly Whyburn if for
every non-closed set $A \subset X$ there is a subset $B \subset A$ such that
$\overline{B} \setminus A$ is a singleton. We prove that every countably
compact Urysohn space of cardinality smaller than the continuum is weakly
Whyburn and show that, consistently, the Urysohn assumption is essential. We
simultaneously solve a question of Pelant, Tkachenko, Tkachuk and Wilson and
one of Bella, Costantini and ourselves by constructing a Lindel\"of $P$-space
of cardinality $\omega_2$ that is not weakly Whyburn. We give conditions for a
weak Whyburn space to be pseudoradial and construct a countably compact weakly
Whyburn non-pseudoradial regular space, which solves a question asked by Bella
in private communication.}

\arXivl{1506.00224}
{A note on local properties in products}
{Paolo Lipparini}
{We give conditions under which a product of topological spaces satisfies some
local property. The conditions are necessary and sufficient when the
corresponding global property is preserved under finite products. Further
examples include local sequential compactness, local Lindel\"ofness, the local
Menger property.}

\arXiv{1506.06080}
{Point-open games and productivity of dense-separable property}
{Jarno Talponen}
{In this note we study the point-open topological games to analyze the least
upper bound for density of dense subsets of a topological space. This way we
may also analyze the behavior of such cardinal invariants in taking products of
spaces. Various related cardinal equalities and inequalities are given. As an
application we take a look at Banach spaces with the property (CSP) which can
be formulated by stating that each weak-star dense linear subspace of the dual
is weak-star separable.}

\arXivl{1507.02134}
{Infinite games and chain conditions}
{Santi Spadaro}
{We apply the theory of infinite two-person games to two well-known problems
in topology: Suslin's Problem and Arhangel'skii's problem on $G_\delta$ covers
of compact spaces. More specifically, we prove results of which the following
two are special cases: 1) every linearly ordered topological space satisfying
the game-theoretic version of the countable chain condition is separable and 2)
in every compact space satisfying the game-theoretic version of the weak
Lindel\"of property, every cover by $G_\delta$ sets has a continuum-sized
subcollection whose union is $G_\delta$-dense.}

\arXiv{1507.02496}
{Nonmeasurable sets and unions with respect to selected ideals especially
  ideals defined by trees}
{Robert Ralowski and Szymon Zeberski}
{In this paper we consider nonmeasurablity with respect to $\sigma$-ideals
defined be trees. First classical example of such ideal is Marczewski ideal
$s_0$. We will consider also ideal $l_0$ defined by Laver trees and $m_0$ defined by
Miller trees. With the mentioned ideals one can consider $s$, $l$ and
$m$-measurablility.

  We have shown that there exists a subset $A$ of the Baire space which is $s$, $l$
and $m$ nonmeasurable at the same time. Moreover, $A$ forms m.a.d.\ family which is
also dominating. We show some examples of subsets of the Baire space which are
measurable in one sense and nonmeasurable in the other meaning.

  We also examine terms nonmeasurable and completely nonmeasurable (with
respect to several ideals with Borel base). There are several papers about
finding (completely) nonmeasurable sets which are the union of some family of
small sets. In this paper we want to focus on the following problem: Let $P$ be
a family of small sets. Is it possible that for all $A$ which is a subset of $P$,
union of $A$ is nonmeasurable implies that union of $A$ is completely
nonmeasurable?
  
  We will consider situations when $P$ is a partition of $R$, $P$ is point-finite
family and $P$ is point-countable family. We give an equivalent statement to CH
using terms nonmeasurable and completely nonmeasurable.}

\arXivl{1507.06717}
{Some observations on the Baireness of $C_k(X)$ for a locally compact space
  $X$}
{Franklin D. Tall}
{We prove some consistency results concerning the Moving Off Property for
locally compact spaces and thus the question of whether their function spaces
are Baire.}

\arXivl{1510.02654}
{Reconstructing Compact Metrizable Spaces}
{Paul Gartside, Max F. Pitz, Rolf Suabedissen}
{The deck, $\mathcal{D}(X)$, of a topological space $X$ is the set
$\mathcal{D}(X)=\{[X \setminus \{x\}]\colon x \in X\}$, where $[Y]$ denotes the
homeomorphism class of $Y$. A space $X$ is (topologically) reconstructible if
whenever $\mathcal{D}(Z)=\mathcal{D}(X)$ then $Z$ is homeomorphic to $X$. It is
known that every (metrizable) continuum is reconstructible, whereas the Cantor
set is non-reconstructible.

  The main result of this paper characterises the non-reconstructible compact
metrizable spaces as precisely those where for each point $x$ there is a
sequence $\langle B_n^x \colon n \in \mathbb{N}\rangle$ of pairwise disjoint
clopen subsets converging to $x$ such that $B_n^x$ and $B_n^y$ are homeomorphic
for each $n$, and all $x$ and $y$.

  In a non-reconstructible compact metrizable space the set of $1$-point
components forms a dense $G_\delta$. For $h$-homogeneous spaces, this condition
is sufficient for non-reconstruction. A wide variety of spaces with a dense
$G_\delta$ set of $1$-point components are presented, some reconstructible and
others not reconstructible.}

\arXivl{1511.07062}
{On topological groups admitting a base at identity indexed with
  $\omega^\omega$}
{Arkady G. Leiderman, Vladimir G. Pestov, and Artur H. Tomita}
{A topological group $G$ is said to have a $\mathfrak G$-base if the
neighbourhood system at identity admits a monotone cofinal map from the
directed set $\omega^\omega$. In particular, every metrizable group is such,
but the class of groups with a $\mathfrak G$-base is significantly wider. The
aim of this article is to better understand the boundaries of this class, by
presenting new examples and counter-examples. Ultraproducts and
non-arichimedean ordered fields lead to natural families of non-metrizable
groups with a $\mathfrak G$-base which nevertheless have the Baire property.
More examples come from such constructions as the free topological group and
the free Abelian topological group of a Tychonoff (more generally uniform)
space, as well as the free product of topological groups. Our results answer
some questions previously stated in the literature.}

%%%%%%%%%%%%%%%%%%%%%%%%%%%%

\section{Short announcements}\label{RA}

\AMS{Strong colorings yield $\kappa$-bounded spaces with discretely
untouchable points}
{Istvan Juhasz; Saharon Shelah}
{http://www.ams.org/journal-getitem?pii=S0002-9939-2014-12394-X}

\arXiv{1501.01949}
{$P$-Paracompact and $P$-Metrizable Spaces}
{Ziqin Feng, Paul Gartside, Jeremiah Morgan}

\arXiv{1501.02877}
{Density character of subgroups of topological groups}
{Arkady Leiderman, Sidney A. Morris, Mikhail G. Tkachenko}

\arXiv{1501.06972}
{Nonseparable growth of the integers supporting a measure}
{Piotr Drygier and Grzegorz Plebanek}

\AMS{Forcing consequences of $PFA$ together with the continuum large}
{David Aspero; Miguel Angel Mota}
{http://www.ams.org/journal-getitem?pii=S0002-9947-2015-06205-9}

\arXiv{1503.04278}
{On the submetrizability number and $i$-weight of quasi-uniform spaces
  and paratopological groups}
{Taras Banakh and Alex Ravsky}

\arXiv{1503.04480}
{Verbal covering properties of topological spaces}
{Taras Banakh and Alex Ravsky}

\AMS{On the collection of Baire class one functions on the irrationals}
{Roman Pol}
{www.ams.org/journal-getitem?pii=S0002-9939-2015-12583-X}

\arXiv{1504.02765}
{Measuring sets with translation invariant Borel measures}
{Andr\'as M\'ath\'e}

\arXiv{1504.01785}
{Cardinalities of weakly Lindel\"of spaces with regular
  $G_\kappa$-diagonals}
{Ivan S. Gotchev}

\arXiv{1504.01790}
{Generalizations of two cardinal inequalities of Hajnal and Juh\'asz}
{Ivan S. Gotchev}

\arXiv{1504.04198}
{Topological properties of function spaces $C_k(X,2)$ over
  zero-dimensional metric spaces $X$}
{S. Gabriyelyan}

\arXiv{1504.04202}
{The Ascoli property for function spaces and the weak topology of Banach
  and Fr\'echet spaces}
{S. Gabriyelyan, J. Kakol, G. Plebanek}

\arXiv{1505.02319}
{Ordering A Square}
{Raushan Z. Buzyakova}

\arXiv{1505.06251}
{Increasing chains and discrete reflection of cardinality}
{Santi Spadaro}

\arXiv{1506.00206}
{Pinning Down versus Density}
{Istv\'an Juh\'asz, Lajos Soukup, Zolt\'an Szentmikl\'ossy}

\arXiv{1506.04665}
{Regular $G_\delta$-diagonals and some upper bounds for cardinality of
  topological spaces}
{Ivan S. Gotchev, Mikhail G. Tkachenko, and Vladimir V. Tkachuk}

\arXiv{1506.08969}
{Cardinal invariants distinguishing permutation groups}
{Taras Banakh and Heike Mildenberger}

\arXiv{1507.06684}
{Cardinality bounds involving the skew-$\lambda$ Lindel\"of degree and
  its variants}
{Nathan Carlson, Jack Porter}

\arXiv{1508.01541}
{Compact spaces with a $\mathbb{P}$-diagonal}
{Alan Dow and Klaas Pieter Hart}

\arXiv{1509.01601}
{Far points and discretely generated spaces}
{Alan Dow, Rodrigo Hern\'andez-Guti\'errez}

\arXiv{1509.02874}
{The weight and Lindel\"of property in spaces and topological groups}
{Mikhail G. Tkachenko}

\arXiv{1509.05282}
{On the problem of compact totally disconnected reflection of
  nonmetrizability}
{Piotr Koszmider}

\arXiv{1509.05542}
{On the sequential closure of the set of continuous functions in the
  space of separately continuous functions}
{Taras Banakh}

\AMS{Classifying invariant $\sigma$-ideals with analytic base on good
Cantor measure spaces}
{Taras Banakh; Robert Ralowski; Szymon Zeberski}
{www.ams.org/journal-getitem?pii=S0002-9939-2015-12709-8}

\section{Problem of the Issue}

The undefined terminology is provided at the end of this section.
The Isbell--Mr\'owka $\Psi$-spaces are classic examples in the realm of topological covering properties.
For obvious reasons (they are not even Lindel\"of), 
Menger's property is not the right notion to consider in the realm of $\Psi$-spaces,
and the correct notion, as observed by Bonanzinga and Matveev~\cite{MilenaMisha},
is the \emph{star-Menger} property.
Initial effort to study the question which $\Psi$-spaces are star-Menger was put forth in
the papers~\cite{MilenaMisha,Ts}. 
Despite these, \emph{it is still unknown whether such spaces consistently exist!}
More precisely, the only examples of star-Menger spaces known thus far
is those of cardinality smaller than $\fd$. Such spaces are star-Menger for an obvious
counting reason, and are thus trivial examples.

\bprb[Bonanzinga--Matveev~\cite{MilenaMisha}]\label{prb:main}
Is there, consistently, a star-Menger $\Psi$-space of cardinality $\ge\fd$?
\eprb

The problem can be formulated in a combinatorial manner, with no referecence to topology.
In particular, the definitions below are not necessary in order to consider this problem.
This is due to the following result.

\bthm[\cite{Ts}]
Let $\cA\sub P(\bbN)$ be an almost disjoint family.
The following assertions are equivalent:
\be
\item The Isbell--Mr\'owka space $\Psi(\cA)$ is star-Menger.
\item For each function $A\mapsto f_A$ from $\cA$ to $\NN$,
there are finite sets $\cF_1,\cF_2,\dots\sub\cA$ such that, for each $A\in\cA$,
there is $n$ with $(A\sm f_A(n))\cap \Un_{B\in\cF_n}(B\sm f_B(n))\neq\emptyset$.
\ee
\ethm

It is only known that there are many cardinal numbers with, provably, no star-Menger $\Psi$-spaces of that cardinality~\cite{MilenaMisha, Ts}. In particular, 
the inequality 
$\fc>\aleph_\omega$ must hold in every model witnessing a positive solution of Problem~\ref{prb:main}.

\subsection*{Basic definitions}
A family $\cA\sub P(\bbN)$ is \emph{almost disjoint} if every element of $\cA$ is infinite, and the sets
$A\cap B$ are finite for all distinct elements $A,B\in \cA$.
For an almost disjoint family $\cA$, let $\Psi(\cA):=\cA\cup\bbN$. A topology on $\Psi(\cA)$ is defined as follows.
The natural numbers are isolated, and for each element $A\in\cA$ and each finite set $F\sub\bbN$,
the set $\{A\}\cup (A\sm F)$ is a basic open neighborhood of $A$. Spaces constructed in this manner
are called \emph{$\Psi$-spaces}.

For a set $X$, a subset $A$ of $X$ and a family $\cU$ of subsets of $X$, let
$\st(A,\cU):=\Un\set{U\in\cU}{A\cap U\ne \emptyset}$.
A topological space $X$ is \emph{star-Lindel\"of} if
every open cover $\cU$ of $X$ has a countable subset $\cV$ such that $X=\st(\Un\cV,\cU)$.
It is easy to see that uncountable $\Psi$-spaces are not Lindel\"of. Being separable, though,
all $\Psi$-spaces are star-Lindel\"of.

\emph{Menger's property} is the following selective version of Lindel\"of's property:
For every sequence $\cU_1,\cU_2,\dots$ of open covers of $X$,
there are finite sets $\cF_1\sub\cU_1, \cF_2\sub\cU_2, \dots$ such that the family 
$\{\Un\cF_1,\Un\cF_2,\dots\}$
covers $X$.

A topological space $X$ is
\emph{star-Menger} if for every sequence $\cU_1,\cU_2,\dots$ of open covers of $X$,
there are finite sets $\cF_1\sub\cU_1, \cF_2\sub\cU_2, \dots$
such that the family $\{\st(\Un\cF_1,\cU_1),\allowbreak\st(\Un\cF_2,\cU_2),\allowbreak\dots\}$
covers $X$.

\ed